\documentclass[11pt]{amsart}
\usepackage{amsmath,amsthm,amsfonts,amssymb,latexsym,epsfig}

\newtheorem{theorem}{Theorem}[section]

\newtheorem{lemma}{Lemma}[section]

\newtheorem{cor}{Corollary}[section]

\numberwithin{equation}{section}

\theoremstyle{definition}

\theoremstyle{remark}

\newcommand{\sumi}{\mathop{{\sum}^{'}}}

\begin{document}
\title{On a result of Levin and Ste\v ckin}
\author{Peng Gao}
\address{Division of Mathematical Sciences, School of Physical and Mathematical Sciences,
Nanyang Technological University, 637371 Singapore}
\email{penggao@ntu.edu.sg}


\begin{abstract}
  We extend a result of Levin and Ste\v ckin concerning an inequality analogous to Hardy's inequality.
\end{abstract}

\maketitle
\section{Introduction}
\label{sec 1} \setcounter{equation}{0}

  Let $p>1$ and $l^p$ be the Banach space of all complex sequences ${\bf a}=(a_n)_{n \geq 1}$. The celebrated
   Hardy's inequality \cite[Theorem 326]{HLP} asserts that for
   $p>1$ and any ${\bf a} \in l^p$,
\begin{equation}
\label{eq:1} \sum^{\infty}_{n=1}\big{|}\frac {1}{n}
\sum^n_{k=1}a_k\big{|}^p \leq \Big(\frac {p}{p-1}
\Big)^p\sum^\infty_{k=1}|a_k|^p.
\end{equation}

   As an analogue of Hardy's inequality,  Theorem
   345 of \cite{HLP} asserts that the following inequality holds for $0<p<1$ and $a_n \geq 0$ with $c_p=p^p$:
\begin{equation}
\label{1}
  \sum^{\infty}_{n=1}\Big( \frac 1{n} \sum^{\infty}_{k=n}a_k \Big
  )^p \geq c_p \sum^{\infty}_{n=1}a^p_n.
\end{equation}
   It is noted in \cite{HLP} that the constant $c_p=p^p$ may not be best possible
  and a better constant was indeed obtained by Levin and Ste\v ckin
  \cite[Theorem 61]{L&S}. Their result is more general as they proved, among other things, the
   following inequality (\cite[Theorem 62]{L&S}), valid for $0<r \leq p
   \leq 1/3$ or $1/3<p<1, r \leq (1-p)^2/(1+p)$ with $a_n \geq 0$,
\begin{equation}
\label{4.1}
  \sum^{\infty}_{n=1}\frac 1{n^r} \Big( \sum^{\infty}_{k=n}a_k \Big
  )^p \geq \Big (\frac {p}{1-r} \Big )^p \sum^{\infty}_{n=1}\frac {a^p_n}{n^{r-p}}.
\end{equation}
   We note here the constant $(p/(1-r))^p$ is best possible, as
   shown in \cite{L&S} by setting $a_n=n^{-1-(1-r)/p-\epsilon}$ and letting $\epsilon \rightarrow
   0^+$. This implies inequality \eqref{1} for $0<p \leq 1/3$ with the
   best possible constant $c_p=(p/(1-p))^p$. On the other hand, it is also easy to see that
  inequality \eqref{1} fails to hold with $c_p=(p/(1-p))^p$
  for $p \geq 1/2$. The point is that in these cases $p/(1-p) \geq 1$ so
  one can easily construct counter examples.

   A simpler proof of Levin and Ste\v ckin's result (for $0< r=p \leq 1/3$) is given in \cite{G6}. It is also pointed out there that using a different approach,
   one may be able to extend
   their result to $p$ slightly larger than $1/3$, an example
   is given for $p=0.34$. The calculation however is more involved and therefore it is desirable to have a simpler approach.
   For this, we let $q$ be the number defined by $1/p+1/q=1$ and note that by the duality principle (see \cite[Lemma 2]{M} but note that our situation is slightly different since
    we have $0<p<1$ with an reversed inequality),
    the case $0<r<1, 0<p<1$ of inequality \eqref{4.1} is equivalent to the
    following one for $a_n>0$:
\begin{equation}
\label{3}
  \sum^{\infty}_{n=1}\Big( n^{(r-p)/p}\sum^{n}_{k=1}\frac {a_k}{k^{r/p}} \Big
  )^{q} \leq  \Big ( \frac {p}{1-r} \Big )^q \sum^{\infty}_{n=1}a^{q}_n.
\end{equation}
    The above inequality can be regarded as an analogue of a result of Knopp \cite{K}, which asserts that Hardy's
    inequality \eqref{eq:1} is still valid for $p<0$ if we assume $a_n>0$.
    We may also regard inequality \eqref{3} as an inequality concerning the factorable matrix with entries $n^{(r-p)/p}k^{-r/p}$ for $k \leq n$ and $0$ otherwise.
    Here we recall that a matrix $A=(a_{nk})$ is factorable if it is a lower triangular matrix
    with $a_{nk}=a_nb_k$ for $1 \leq k \leq n$. We note that the approach in \cite{G7} for the $l^p$ norms of weighted mean
    matrices can also be easily adopted to treat the $l^p$ norms
    of factorable matrices and it is our goal in this paper to use
   this similar approach to extend the result of Levin and Ste\v
   ckin. Our main result is
\begin{theorem}
\label{thm6}
  Inequality \eqref{1}
  holds with the best possible constant $c_p=(p/(1-p))^p$ for any $1/3<p<1/2$ satisfying
\begin{equation}
\label{1.4}
   2^{p/(1-p)}\Big (\Big(\frac {1-p}{p} \Big
)^{1/(1-p)}-\frac {1-p}{p}\Big ) -(1+\frac {3-1/p}2)^{1/(1-p)}
\geq 0.
\end{equation}
   In particular, inequality \eqref{1}
  holds for $0 < p  \leq 0.346$.
\end{theorem}

     It readily follows from Theorem \ref{thm6} and our discussions above that we have the
     following dual version of Theorem \ref{thm6}:
\begin{cor}
\label{dual}
   Inequality \eqref{3} holds with $r=p$ for any $1/3<p<1/2$ satisfying \eqref{1.4} and the constant is
   best possible. In particular, inequality \eqref{3} holds with
   $r=p$ for $0 < p  \leq 0.346$.
\end{cor}

     An alternative proof of Theorem \ref{thm6} is given in Section \ref{sec 3},
     via an approach using the duality principle. In Section \ref{sec4}, we shall study some inequalities
     which can be regarded as generalizations of \eqref{1}.
     Motivations for considerations for such inequalities
     come both from their integral analogues as well as from their
     counterparts in the $l^p$ spaces. As an example, we consider the following inequality for $0<p<1, 0<\alpha<1/p$:
\begin{equation}
\label{1.5}
   \sum^{\infty}_{n=1} \Big(
\sum^{\infty}_{k=n}\frac {\alpha k^{\alpha-1}}{n^{\alpha}}a_k \Big
  )^p \geq \Big ( \frac {\alpha p}{1-\alpha p} \Big )^p \sum^{\infty}_{n=1}a^p_n.
\end{equation}
    As in the case of \eqref{1}, the above inequality doesn't hold for all $0<p<1,
    0<\alpha<1/p$. In Section \ref{sec4}, we will however prove a result concerning the validity of \eqref{1.5} that can be regarded
     as an analogue to
    that of Levin and Ste\v ckin's concerning the validity of
    \eqref{1}.

    Inequality \eqref{1.5} is motivated partially by integral
    analogues of \eqref{1}, as we shall explain in Section
    \ref{sec4}. It is also motivated by the following inequality
    for $p>1, \alpha p>1$, $a_n \geq 0$:
\begin{equation}
\label{1.2}
  \sum^{\infty}_{n=1}\Big ( \sum^n_{k=1}\frac {\alpha
  k^{\alpha-1}}{n^{\alpha}}a_k \Big )^p \leq \Big (\frac
  {\alpha p}{\alpha p-1} \Big )^p\sum^{\infty}_{n=1}a^p_n.
\end{equation}
    The above inequality is in turn motivated by the following inequality
\begin{equation}
\label{1.3}
  \sum^{\infty}_{n=1}\Big ( \sum^n_{k=1}\frac {
  k^{\alpha-1}}{\sum^n_{i=1}i^{\alpha-1}}a_k \Big )^p \leq \Big (\frac
  {\alpha p}{\alpha p-1} \Big )^p\sum^{\infty}_{n=1}a^p_n.
\end{equation}
   Inequality \eqref{1.3} was first suggested by Bennett \cite[p.
   40-41]{B4}, see \cite{G8} and the references therein for recent
   progress on this. We point out here that it is easy to see that inequality
   \eqref{1.2} implies \eqref{1.3} when $\alpha>1$, hence it is
   interesting to know for what $\alpha$'s, inequality \eqref{1.2} is valid. We first note that on setting $a_1=1$ and $a_n=0$
   for $n \geq 2$ in \eqref{1.2} that it is impossible for it to
   hold when $\alpha$ is large for fixed $p$. On the other hand,
   when $\alpha=1$, inequality \eqref{1.2} becomes Hardy's
   inequality and hence one may expect it to hold for $\alpha$
   close to $1$ and we shall
   establish such a result in Section
   \ref{sec3}.
\section{Proof of Theorem \ref{thm6}}
\label{sec 2'} \setcounter{equation}{0}

   First we need a lemma:
\begin{lemma}
\label{lem1}
   The following inequality holds for $0
   \leq y \leq 1$ and $1/2<t<1$:
\begin{equation}
\label{2.1}
   (1+y/(2t))^{1+t}-(1+y)^{-t}(1+(2t-1)y/(2t))^{1+t}- y/t  \geq 0.
\end{equation}
\end{lemma}
\begin{proof}
   We set $x=y/(2t)$ so that $0 \leq x \leq 1$ and we recast
   the above inequality as
\begin{equation*}
   f(x,t):=(1+x)^{1+t}-(1+2tx)^{-t}(1+(2t-1)x)^{1+t}- 2x  \geq 0.
\end{equation*}
   Direct calculation shows that $f(0,t)=\frac {\partial f}{\partial x} (0,t)=0$ and
\begin{equation*}
   \frac {\partial^2 f}{\partial x^2} (x,t)=t(1+t)(1+x)^{t-1}\Big (1-(1+2tx)^{-t-2}(1+(2t-1)x)^{t-1}(1+x)^{1-t} \Big ) :=t(1+t)(1+x)^{t-1}g(x,t).
\end{equation*}
   Note that
\begin{equation*}
   \frac {\partial g}{\partial x } (x,t)=(1+2tx)^{-t-3}(1+(2t-1)x)^{t-2}(1+x)^{-t}\Big (2(4t-1)+4t(4t-1)x+2t(2t-1)(t+2)x^2) \geq 0 .
\end{equation*}
  As $g(0, t)=0$, it follows that $g(x, t) \geq 0$ for $0 \leq x
  \leq 1$ which in turn implies the assertion of the lemma.
\end{proof}

     We now describe a general approach towards establishing inequality \eqref{4.1} for $0<r<1, 0 < p < 1$. A modification from the approach in Section 3
of \cite{G6} shows
   that in order for \eqref{4.1} to hold for any given $p$ with
   $c_{p,r} (=(p/(1-r))^p)$, it suffices to find
  a sequence ${\bf w}$ of positive terms for each $0<r<1$ and $0<p < 1$, such
  that for any integer $n \geq 1$,
\begin{equation*}
 n^{(p-r)/(1-p)}(w_1+\cdots+w_n)^{-1/(1-p)} \leq
 c^{-1/(1-p)}_{p,r}
 \Big( \frac {w^{-1/(1-p)}_n}{n^{r/(1-p)}}- \frac
 {w^{-1/(1-p)}_{n+1}}{(n+1)^{r/(1-p)}}\Big).
\end{equation*}
    We note here that if we study the equivalent inequality \eqref{3} instead, then we can also obtain the above inequality from inequality
    (2.2) of \cite{G6}, on setting $\Lambda_n=n^{-(r-p)/p}, \lambda_n=n^{-r/p}$ there. For the moment, we assume $c_{p,r}$ is an
arbitrary fixed positive number and on setting
$b^{p-1}_n=w_n/w_{n+1}$, we can recast the above inequality as
\begin{equation*}
  \Big(\sum^n_{k=1}\prod^n_{i=k}b^{p-1}_i \Big )^{-1/(1-p)} \leq
 c^{-1/(1-p)}_{p,r}n^{(r-p)/(1-p)}
 \Big( \frac {b_n}{n^{r/(1-p)}}- \frac
 {1}{(n+1)^{r/(1-p)}}\Big).
\end{equation*}
  The choice of $b_n$ in Section 3 of \cite{G6} suggests that for optimal choices of the $b_n$'s, we
  may have asymptotically $b_n \sim 1+c/n$ as $n \rightarrow +\infty$ for
  some positive constant $c$ (depending on $p$). This observation implies that
  $n^{1/(1-p)}$ times the right-hand side expression above should
  be asymptotically a constant. To take the advantage of possible
  contributions of higher order terms, we now further recast the above inequality as
\begin{equation}
\label{4.2}
 \Big(\frac {1}{n+a}\sum^n_{k=1}\prod^n_{i=k}b^{p-1}_i \Big )^{-1/(1-p)} \leq
 c^{-1/(1-p)}_{p,r}n^{(r-p)/(1-p)}(n+a)^{1/(1-p)}
 \Big( \frac {b_n}{n^{r/(1-p)}}- \frac
 {1}{(n+1)^{r/(1-p)}}\Big),
\end{equation}
  where $a$ is a constant (may depend on $p$) to be chosen later. It will also be clear from our arguments below that the choice of $a$ will not affect
  the asymptotic behavior of $b_n$ to the first order of magnitude. We now choose $b_n$ so that
\begin{equation}
\label{2.3}
 n^{(r-p)/(1-p)}(n+a)^{1/(1-p)}
 \Big( \frac {b_n}{n^{r/(1-p)}}- \frac
 {1}{(n+1)^{r/(1-p)}}\Big)=c^{-\alpha/(1-p)}_{p,r},
\end{equation}
  where $\alpha$ is a parameter to be chosen later.
  This implies that
\begin{equation*}
   b_n=c^{-\alpha/(1-p)}_{p,r}\frac {n^{p/(1-p)}}{(n+a)^{1/(1-p)}}+\frac {n^{r/(1-p)}}{(n+1)^{r/(1-p)}}.
\end{equation*}
   For the so chosen $b_n$'s, inequality \eqref{4.2} becomes
\begin{equation}
\label{4.3}
  \sum^n_{k=1}\prod^n_{i=k}b^{p-1}_i  \geq (n+a)c^{1+\alpha}_{p,r}.
\end{equation}
   We first assume the above inequality holds for $n=1$. Then
   induction shows it holds for all $n$ as long as
\begin{equation*}
   b^{1-p}_n  \leq
 \frac {n+a+c^{-(1+\alpha)}_{p,r}-1}{n+a}.
\end{equation*}
  Taking account into the value of $b_n$, the above becomes (for $0 \leq y \leq 1$ with
  $y=1/n$)
\begin{equation}
\label{4.4}
   (1+(a+c^{-(1+\alpha)}_{p,r}-1)y)^{1/(1-p)}-(1+y)^{-r/(1-p)}(1+ay)^{1/(1-p)}- c^{-\alpha/(1-p)}_{p,r}y  \geq 0.
\end{equation}
   The first order term of the Taylor expansion of
   the left-hand side expression above implies that it is necessary to have
\begin{equation*}
    c^{-(1+\alpha)}_{p,r}-(1-p)c^{-\alpha/(1-p)}_{p,r}+r-1 \geq 0.
\end{equation*}
    For fixed $c_{p,r}$, the left-hand side expression above is
    maximized when $\alpha=1/p-1$ with value $pc^{-1/p}_{p,r}+r-1$.
    This suggests us to take $c_{p,r}=(p/(1-r))^p$. From now on we
    fix $c_{p,r}=(p/(1-r))^p$ and note that in this case \eqref{4.4}
    becomes
\begin{equation}
\label{4.5}
   (1+(a+(1-r)/p-1)y)^{1/(1-p)}-(1+y)^{-r/(1-p)}(1+ay)^{1/(1-p)}- \frac {1-r}{p}y  \geq 0.
\end{equation}

    We note that the choice of $a=0$ in \eqref{4.5} with $r=p$ reduces to that considered in Section 3 of
  \cite{G6} (in which case the case $n=1$ of \eqref{4.3} is also included in \eqref{4.5}).
  Moreover, with $a=0$ in the above inequality and following the treatment in Section 3 of \cite{G6}, one is
    able to improve some cases of
    the above mentioned result of Levin and Ste\v ckin concerning the validity of \eqref{4.1}.
    We shall postpone the discussion of this to the next section and focus now on the proof of Theorem \ref{thm6}.
  Since the cases $0< p \leq 1/3$ of the assertion of the theorem
  are known, we may assume $1/3<p<1/2$ from now on. In this case we set $r=p$ in \eqref{4.5} and
  Taylor expansion shows that it is necessary to have $a \geq (3-1/p)/2$ in order
  for
  inequality \eqref{4.5} to hold.  We now take $a=(3-1/p)/2$ and write $t=p/(1-p)$ to
  see that inequality \eqref{4.5} is reduced to \eqref{2.1} and
  Lemma \ref{lem1} now implies that inequality \eqref{4.5} holds
  in this case. Inequality \eqref{1} with the best possible constant $c_{p}=(p/(1-p))^p$
  thus follows for any $1/3<p<1/2$ as long as the case $n=1$ of \eqref{4.3} is satisfied, which is just inequality \eqref{1.4} and
  this proves the first assertion of Theorem \ref{thm6}. For the second assertion, we note that inequality \eqref{1.4}
  can be rewritten as
\begin{equation}
\label{2.7}
   \frac {2^t}{t}(t^{-t}-1) \geq (1+a)^{1/(1-p)},
\end{equation}
   where $t$ is defined as above. Note that $1/2<t<1$ for $1/3 < p <1/2$ and
   both $2^t/t$ and $t^{-t}-1$ are decreasing functions of $t$. It
   follows that the left-hand side expression of \eqref{2.7} is a
   decreasing function of $p$.
   Note also that for fixed $a$, the right-hand side expression of \eqref{2.7} is an increasing
   function of $p<1$. As $a=(3-1/p)/2$ in our case, it follows that one just need to check the above inequality for $p=0.346$ and the
   assertion of the theorem now follows easily.

   We remark here that in the proof of Theorem \ref{thm6}, if instead of choosing $b_n$ to satisfy
   \eqref{2.3} (with $r=p$ and $c_{p,p}=(p/(1-p))^p$ there), we
   choose $b_n$ for $n \geq 2$ so that
\begin{equation*}
   (n+c)^{1/(1-p)}
 \Big( \frac {1}{n^{p/(1-p)}}- \frac
 {1}{(n+1)^{p/(1-p)}b_n}\Big)=(1-p)/p.
\end{equation*}
    Moreover, note that we can also rewrite \eqref{4.2} for $n \geq 2$ as (with
    $a$ replaced by $c$ and $r=p$, $c_{p,p}=(p/(1-p))^p$)
\begin{equation*}
 \Big(\frac {1}{n+c}\Big( \sum^{n-1}_{k=1}\prod^{n-1}_{i=k}b^{p-1}_i+1 \Big )\Big )^{-1/(1-p)} \leq
 \Big (\frac {1-p}{p} \Big )^{p/(1-p)}(n+c)^{1/(1-p)}
 \Big( \frac {1}{n^{p/(1-p)}}- \frac
 {1}{(n+1)^{p/(1-p)}b_n}\Big).
\end{equation*}
  If we further choose $b_1$ so that
\begin{equation*}
 1=\Big (\frac {1-p}{p} \Big )^{p/(1-p)}
 \Big( 1- \frac
 {1}{2^{p/(1-p)}b_1}\Big).
\end{equation*}
   Then repeating the same process as in the proof of Theorem
   \ref{thm6}, we find that the induction part (with $c=(1/p-1)/2$ here)
   leads back to inequality \eqref{4.5} (with $r=p$ and $a=(3-1/p)/2$ there) while the initial
   case (corresponding to
   $n=2$ here) is just \eqref{2.7}, so this approach gives another
   proof of Theorem \ref{thm6}.

   We end this section by pointing out the relation between the treatment in
   Sections
   3 and 4 in \cite{G6} on inequality \eqref{1}. We note that it is shown in Section 3 of \cite{G6}
   that for any $N \geq 1$ and any positive sequence ${\bf w}$, we
   have
\begin{equation}
\label{5.1}
   \sum^{N}_{n=1}a^p_n \leq \sum^{N}_{n=1}w_n\Big ( \sum^{N}_{k=n}\Big(\sum^k_{i=1}w_i \Big)^{-1/(1-p)}\Big
  )^{1-p} \Big ( \sum^{N}_{k=n}a_k\Big
  )^{p}.
\end{equation}
    We now use $w_n=\sum^n_{k=1}w_k-\sum^{n-1}_{k=1}w_k$ and set
    (with $\nu_{N}=0$)
\begin{equation*}
     \nu_n=\frac {\sum^{N}_{k=n+1}\Big(\sum^k_{i=1}w_i
    \Big)^{1/(p-1)}}{\Big(\sum^n_{i=1}w_i
    \Big)^{1/(p-1)}}
\end{equation*}
    to see that inequality \eqref{5.1} leads to (with $\nu_0=0$)
\begin{equation*}
   \sum^{N}_{n=1}a^p_n \leq \sum_{n=1}^{N}\Big ( (1+\nu_n)^{1-p}-
   \nu^{1-p}_{n-1} \Big )\Big ( \sum^{N}_{k=n}a_k\Big
  )^{p}.
\end{equation*}
   The above inequality is essentially what's used in Section 4 of
   \cite{G6}.

\section{An Alternative Proof of Theorem \ref{thm6}}
\label{sec 3} \setcounter{equation}{0}
    In this section we give an alternative proof of Theorem \ref{thm6},
    using the following:
\begin{lemma}[{\cite[Lemma 2.4]{G}}]
\label{lem6.1}
  Let $\{ \lambda_i \}^{\infty}_{i \geq 1}, \{ a_i \}^{\infty}_{i \geq 1}$ be two sequences of positive real numbers and
  let $S_n=\sum_{i=1}^n \lambda_i
  a_i$. Let $0 \neq p<1$ be fixed and let $\{ \mu_i \}^{\infty}_{i \geq 1}, \{ \eta_i \}^{\infty}_{i \geq 1}$ be two
  positive sequences of real numbers such
  that $\mu_i \leq \eta_i$ for $0<p<1$ and $\mu_i \geq \eta_i$ for
  $p<0$, then for $n \geq 2$,
\begin{equation}
\label{6.1}
   \sum_{i=2}^{n-1}\Big ( \mu_i-(\mu^q_{i+1}-\eta^q_{i+1})^{1/q} \Big )S_i^{1/p}+\mu_nS_n^{1/p}
   \leq (\mu^q_{2}-\eta^q_{2})^{1/q}\lambda^{1/p}_1a_1^{1/p}
   +\sum_{i=2}^n \eta_i \lambda^{1/p}_i a_i^{1/p}.
\end{equation}
\end{lemma}

   Following the treatment in Section 4 of \cite{G6}, on first setting $\eta_i =
   \lambda^{-1/p}_i$, then a change of variables: $\mu_i \mapsto
   \mu_i\eta_i$ followed by setting $\mu^q_{i}-1=\nu_i$
   and lastly a further change of variable: $p \mapsto 1/p$, we
   can transform inequality \eqref{6.1} to the following inequality (with $\nu_1=0$ here):
\begin{equation*}
\label{6.2}
    \sum_{i=1}^{n-1}\Big ( \frac {(1+\nu_i)^{1-p}}{\lambda^{p}_i}-
   \frac {\nu^{1-p}_{i+1}}{\lambda^{p}_{i+1}} \Big )S_i^{p}+ \frac {(1+\nu_n)^{1-p}}{\lambda^{p}_n}S_n^{p}
   \leq \sum_{i=1}^n  a_i^{p}.
\end{equation*}
    Here the $\nu_i$'s are arbitrary non-negative real numbers for $2 \leq i \leq n$. On setting
    $\nu_{n+1}$ to be any non-negative real number, we deduce immediately from the above inequality the following:
\begin{equation}
\label{3.02}
    \sum_{i=1}^{n}\Big ( \frac {(1+\nu_i)^{1-p}}{\lambda^{p}_i}-
   \frac {\nu^{1-p}_{i+1}}{\lambda^{p}_{i+1}} \Big )S_i^{p} \leq \sum_{i=1}^n  a_i^{p}.
\end{equation}

    Now we consider establishing inequality \eqref{4.1} for $0<r<1, 1/3<p<1/2$ in general
    and as has been pointed out in Section \ref{sec 1}, we know
    this is equivalent to establishing inequality \eqref{3}.
    Now in order to establish inequality \eqref{3}, it suffices
    to consider the cases of \eqref{3} with the infinite summations replaced by any finite summations,
    say from $1$ to $N \geq 1$ there.
    We now set $n=N, p=q, \lambda_i=i^{-r/p}$ in inequality \eqref{3.02}
    to recast it as (with $\nu_1=0$, $S_n=\sum^{n}_{k=1}k^{-r/p}a_k$ here):
\begin{equation*}
    \sum_{n=1}^{N}\Big ( \frac {(1+\nu_n)^{1/(1-p)}}{n^{r/(1-p)}}-
   \frac {\nu^{1/(1-p)}_{n+1}}{(n+1)^{r/(1-p)}} \Big )S_n^{q} \leq \sum_{n=1}^N  a_n^{q}.
\end{equation*}
     Comparing the above inequality with \eqref{3}, we see that inequality \eqref{3} holds as long as we can find
     non-negative $\nu_n$'s (with $\nu_1=0$) such that
\begin{equation}
\label{3.03}
    \frac {(1+\nu_n)^{1/(1-p)}}{n^{r/(1-p)}}-
   \frac {\nu^{1/(1-p)}_{n+1}}{(n+1)^{r/(1-p)}} \geq n^{(p-r)/(1-p)}\left (\frac {p}{1-r}
   \right )^{p/(1-p)}.
\end{equation}
    Now, on setting for $n \geq 2$,
\begin{equation*}
    \nu_n=\frac {n+a-1}{(1-r)/p},
\end{equation*}
    and $y=1/n$, we see easily that inequality \eqref{3.03} can be
    transformed into inequality \eqref{4.5}. In the case of $r=p$, we further set $a=(3-1/p)/2$ to
    see that the validity of \eqref{4.5} established for this case in Section \ref{sec 2'} ensures the validity
    of \eqref{3.03} for $n \geq 2$. Moreover, with the above chosen
    $\nu_2$ with $r=p$ and $a=(3-1/p)/2$,
    the $n=1$ case of \eqref{3.03} is easily seen to be equivalent
    to inequality \eqref{1.4} and hence this provides an alternative proof of
    Theorem \ref{thm6}.

\section{A generalization of Theorem \ref{thm6}}
\label{sec4} \setcounter{equation}{0}

   For $0<p<1$, let $f(x) \geq 0$ and $\alpha$ be a real number such that $\alpha <1/p$, we have the following
   identity:
\begin{equation}
\label{20}
   \int^{\infty}_{0}\Big (\frac
   {1}{x^{\alpha}}\int^{\infty}_xf(t)t^{\alpha-1}dt
   \Big)^pdx=\Big (\frac {p}{1-\alpha p} \Big )\int^{\infty}_0\Big
   (\frac {1}{x^{\alpha}}\int^{\infty}_{x}f(t)t^{\alpha-1}dt\Big
   )^{p-1}f(x)dx.
\end{equation}
   In the above expression, we assume $f$ is taken so that all the integrals converge. The case of $\alpha=1$ is given in the proof of Theorem
   337 of \cite{HLP} and the general case is obtained by some
   changes of variables. As in the proof of Theorem 337 of
   \cite{HLP}, we then deduce the following inequality (with the same assumptions as above):
\begin{equation*}
   \int^{\infty}_{0}\Big (\frac
   {1}{x^{\alpha}}\int^{\infty}_xf(t)t^{\alpha-1}dt
   \Big)^pdx \geq \Big (\frac {p}{1-\alpha p} \Big )^p\int^{\infty}_0f^p(x)dx.
\end{equation*}
  The above inequality can also be deduced from Theorem 347 of
  \cite{HLP} (see also \cite[(2.4)]{G5}).
   Following the way how Theorem 338 is deduced from Theorem 337 of \cite{HLP}, we deduce similarly from \eqref{20} the following
   inequality for $0<p<1, 0< \alpha <1/p$ and $a_n \geq 0$:
\begin{equation*}
   \sumi^{\infty}_{n=1} \Big( \frac 1{n^{\alpha}} \sum^{\infty}_{k=n}\Big((k+1)^{\alpha}-k^{\alpha}\Big )a_k \Big
  )^p \geq \Big ( \frac {\alpha p}{1-\alpha p} \Big )
  \sum^{\infty}_{n=1}\Big( \frac 1{n^{\alpha}} \sum^{\infty}_{k=n}\Big((k+1)^{\alpha}-k^{\alpha}\Big )a_k \Big
  )^{p-1}a_n.
\end{equation*}
  The dash over the summation on the left-hand side expression above (and in
what follows) means that
   the term for which $n=1$ is to be multiplied by $1+1/(1-\alpha p)$ and
   the constant is best possible (on taking $a_n=n^{-1/p-\epsilon}$ and letting $\epsilon \rightarrow 0^+$).
   The above inequality readily implies the following one by H\"older's
   inequality:
\begin{equation*}
\   \sumi^{\infty}_{n=1} \Big( \frac 1{n^{\alpha}}
\sum^{\infty}_{k=n}\Big((k+1)^{\alpha}-k^{\alpha}\Big )a_k \Big
  )^p \geq \Big ( \frac {\alpha p}{1-\alpha p} \Big )^p \sum^{\infty}_{n=1}a^p_n.
\end{equation*}

    We are thus motivated to consider the above inequality with the dash sign removed and this can be regarded as an analogue of
    inequality \eqref{1} with $c_p=(p/(1-p))^p$, which corresponds to the case $\alpha=1$ here. As in the case of
\eqref{1}, such an inequality
   does not hold for all $\alpha$ and $p$ satisfying $0<p<1$ and
   $0< \alpha <1/p$. However, on setting $a_n=n^{-1/p-\epsilon}$ and
   letting $\epsilon \rightarrow 0^+$, one sees easily that if
   such an inequality holds for
   certain $\alpha$ and $p$, then the constant is best possible. More generally, we can consider the following inequality:
\begin{equation}
\label{21}
 \sum^{\infty}_{n=1}\Big ( \frac
1{\sum^n_{i=1}L^{\alpha-1}_{\beta}(i,
i-1)}\sum^{\infty}_{k=n}L^{\alpha-1}_{\beta}(k \pm 1, k)a_k \Big
)^p \geq \Big ( \frac {\alpha p}{1-\alpha p} \Big
)^p\sum^{\infty}_{n=1}a_n^p,
\end{equation}
   where the function $L_r(a,b)$ for $a>0, b>0, a \neq b$ and $r \neq 0, 1$ (the only case we shall concern here) is defined
   as $L^{r-1}_r(a,b)=(a^r-b^r)/(r(a-b))$. It is known \cite[Lemma 2.1]{alz1.5} that
   the function $r \mapsto L_r(a,b)$ is strictly increasing on ${\mathbb
   R}$. Here we restrict our attention to the plus sign in \eqref{21} for the case
   $\beta>0, \max (1, \beta) \leq \alpha$ and to the minus sign in \eqref{21} for the case $0<
   \alpha <1$ and $\beta \geq \alpha$. Our remark above implies
   that in either case (note that $L_{\beta}(1,0)$ is meaningful)
\begin{equation*}
  \sum^n_{i=1}L^{\alpha-1}_{\beta}(i,
i-1) \leq  \sum^n_{i=1}L^{\alpha-1}_{\alpha}(i,
i-1)=n^{\alpha}/\alpha.
\end{equation*}
   As we also have $L^{\alpha-1}_{\beta}(k \pm 1, k) \geq
   k^{\alpha-1}$, we see that the validity of \eqref{21} follows
   from that of \eqref{1.5}. We therefore focus on \eqref{1.5}
   from now on and we proceed as in Section 3 of \cite{G6} to see that in order
   for inequality \eqref{1.5} to hold, it suffices to find
  a sequence ${\bf w}$ of positive terms for each $0<p < 1$, such
  that for any integer $n \geq 1$,
\begin{equation}
\label{35}
  \Big(\sum^n_{k=1}w_k \Big )^{1/(p-1)} \leq
 \Big (\frac {\alpha p}{1-\alpha p} \Big )^{p/(p-1)}\Big (\alpha n^{\alpha-1} \Big )^{p/(1-p)}
 \Big( \frac {w^{1/(p-1)}_n}{n^{\alpha p/(1-p)}}- \frac
 {w^{1/(p-1)}_{n+1}}{(n+1)^{\alpha p/(1-p)}}\Big).
\end{equation}
   We now choose ${\bf w}$ inductively by setting $w_1=1$ and for
   $n \geq 1$,
\begin{equation*}
   w_{n+1}=\frac {n+1/p-\alpha-1}{n}w_n.
\end{equation*}
   The above relation implies that
\begin{equation*}
  \sum^n_{k=1}w_k=\frac {n+1/p-\alpha-1}{1/p-\alpha}w_n.
\end{equation*}
    We now assume $0<p<1/2$ and note that for the so chosen ${\bf w}$, inequality \eqref{35} follows (with
    $x=1/n$) from $f(x) \geq 0$ for $0 \leq x \leq 1$, where
\begin{equation}
\label{3.5'}
  f(x)=\Big(1+(1/p-\alpha-1)x \Big )^{1/(1-p)}-(1+x)^{-\alpha
  p/(1-p)}-\frac {1-\alpha p}{p}x.
\end{equation}
   As $f(0)=f'(0)=0$, it suffices to show $f''(x) \geq 0$, which
   is equivalent to showing $g(x) \geq 0$ where
\begin{equation*}
  g(x)=\Big(\frac {(1/p-\alpha-1)^2}{\alpha ((\alpha-1)p+1)} \Big
  )^{(1-p)/(1-2p)}(1+x)^{(2+(\alpha-2)p)/(1-2p)}-(1+(1/p-\alpha-1)x).
\end{equation*}
   Now
\begin{align}
\label{36}
  g'(x) &=\Big(\frac {(1/p-\alpha-1)^2}{\alpha ((\alpha-1)p+1)} \Big
  )^{(1-p)/(1-2p)}\Big (\frac {2+(\alpha-2)p}{1-2p} \Big
  )(1+x)^{(2+(\alpha-2)p)/(1-2p)-1}-(1/p-\alpha-1)\\
   & \geq \Big(\frac {(1/p-\alpha-1)^2}{\alpha ((\alpha-1)p+1)} \Big
  )^{(1-p)/(1-2p)}\Big (\frac {2+(\alpha-2)p}{1-2p} \Big
  )-(1/p-\alpha-1):=h(\alpha, p). \nonumber
\end{align}
   Suppose now $\alpha \geq 1$, then when $1/p \geq (\alpha+2)(\alpha+1)/2$, we have $1/p \geq \alpha(\alpha-1)p+2\alpha+1$ since $p < 1/2$ so that
   both inequalities $1/p-\alpha-1 \geq 1$ and $1/p-\alpha-1 \geq \alpha
   ((\alpha-1)p+1)$ are satisfied. In this case we have
\begin{align*}
   h(\alpha, p) & \geq \Big(\frac {(1/p-\alpha-1)^2}{\alpha ((\alpha-1)p+1)} \Big
  )^{(1-p)/(1-2p)}-(1/p-\alpha-1) \geq  \frac {(1/p-\alpha-1)^2}{\alpha ((\alpha-1)p+1)} -(1/p-\alpha-1) \geq 0.
\end{align*}
   It follows $g'(x) \geq 0$ and as $g(0) \geq 0$, we conclude that $g(x) \geq 0$ and hence $f(x) \geq
   0$. Similar discussion leads to the same conclusion for $0< \alpha< 1$ when $p \leq 1/(\alpha+2)$. We now summarize our discussions above in the following
\begin{theorem}
   Let $0<p<1/2$ and $0 < \alpha < 1/p$. Let $h(\alpha, p)$ be defined
   as in \eqref{36}. Inequality \eqref{1.5} holds for $\alpha, p$ satisfying $h(\alpha, p) \geq 0$. In
   particular, when $\alpha \geq 1$,
   inequality \eqref{1.5} holds for $0<p \leq 2/((\alpha+2)(\alpha+1))$. When $0< \alpha \leq 1$,
   inequality \eqref{1.5} holds for $0<p \leq 1/(\alpha+2)$.
\end{theorem}

\begin{cor}
   Let $0<p<1/2$ and $0<  \alpha < 1/p$. Let $h(\alpha, p)$ be defined
   as in \eqref{36}. When $\beta>0, \max (1 , \beta) \leq \alpha$, inequality \eqref{21} holds (where we take the plus sign) for $\alpha, p$ satisfying $h(\alpha, p) \geq 0$. In
   particular,
   inequality \eqref{21} holds for $0<p \leq
   2/((\alpha+2)(\alpha+1))$. When $0<\alpha<1, \beta \geq \alpha$, inequality \eqref{21} holds (where we take the minus sign) for $\alpha, p$ satisfying $h(\alpha, p) \geq 0$. In
   particular,
   inequality \eqref{21} holds for $0<p \leq 1/(\alpha+2)$.
\end{cor}

    We note here a special case of the above corollary, the case $0<\alpha<1$ and $\beta \rightarrow +\infty$  leads to the
    following inequality, valid for $0< p \leq 1/(\alpha+2)$:
\begin{equation*}
   \sum^{\infty}_{n=1} \Big(\frac 1{ \sum^{n}_{i=1}i^{\alpha-1}}
\sum^{\infty}_{k=n}k^{\alpha-1}a_k \Big
  )^p \geq \Big ( \frac {\alpha p}{1-\alpha p} \Big )^p
  \sum^{\infty}_{n=1}a^p_n.
\end{equation*}

We further note here if we set $r=\alpha p$ and $a=0$ in
inequality \eqref{4.5}, then it is reduced to $f(x) \geq 0$ for
$f(x)$
    defined as in \eqref{3.5'}. Since the case $0<r<p \leq 1/3$
    is known, we need only concern the case $\alpha \geq 1$ here
    and we now have the following improvement of the result of Levin
    and Ste\v ckin \cite[Theorem 62]{L&S}:
\begin{cor}
   Let $0<p<1/2$ and $1 \leq  \alpha < 1/p$. Let $h(\alpha, p)$ be defined
   as in \eqref{36}. Inequality \eqref{4.1} holds for $r=\alpha p$
   for $\alpha, p$ satisfying $h(\alpha, p) \geq 0$. In
   particular, inequality \eqref{4.1} holds for $r=\alpha p$
   for $\alpha, p$ satisfying $0<p \leq
   2/((\alpha+2)(\alpha+1))$.
\end{cor}

    Just as Theorem \ref{thm6} and Corollary \ref{dual} are dual versions to each other, our results above can also be stated
    in terms of their dual versions and we shall leave the formulation of the corresponding
    ones to the reader.
\section{Some results on $l^p$ norms of factorable matrices }
\label{sec3} \setcounter{equation}{0}
  In this section we first state some results concerning the $l^p$ norms
    of factorable matrices. In order to compare our result to that of weighted mean matrices, we consider the following type of inequalities:
\begin{equation}
\label{3.1}
   \sum^{\infty}_{n=1}\Big (\sum^n_{k=1}\frac
   {\lambda_k}{\Lambda_n}a_k \Big )^p \leq
   U_p\sum^{\infty}_{n=1}a^p_n,
\end{equation}
   where $p>1, U_p$ is a constant depending on $p$.
    Here we assume the two positive sequences $(\lambda_n)$ and
    $(\Lambda_n)$ are independent (in particular, unlike in the
    weighted mean matrices case, we do not have
    $\Lambda_n=\sum^n_{k=1}\lambda_k$ in general).
    We begin with the following result concerning the bound for
    $U_p$:
\begin{theorem}
\label{thm3.1}
  Let $1<p<\infty$ be fixed in \eqref{3.1}. Let $a$ be a constant such that $\Lambda_n+a\lambda_n>0$ for all $n \geq 1$.
  Let $0<L<p$ be a positive constant and let
\begin{equation*}
  b_n=( \frac {p-L}{p} )(1+a\frac
  {\lambda_n}{\Lambda_n})^{p-1}\frac {\lambda_n}{\Lambda_n}+\frac
  {\lambda_n}{\lambda_{n+1}}.
\end{equation*}
  If for any integer $n \geq 1$, we have
\begin{equation*}
    \sum^n_{k=1}\lambda_k\prod^{n}_{i=k}b^{1/(p-1)}_i \leq
     \frac {p}{p-L}(\Lambda_n+ a \lambda_n),
\end{equation*}
    then inequality \eqref{3.1} holds with $U_{p} \leq (p/(p-L))^p$.
\end{theorem}

   We point out that the proof of the above theorem is analogue to
   that of Theorem 3.1 of \cite{G7}, except instead of choosing
   $b_n$ to satisfy the equation (3.4) in \cite{G7}, we choose
   $b_n$ so that
\begin{equation*}
     \Big( \frac {b_n}{\lambda_n}-\frac {1}{\lambda_{n+1}} \Big )
    \Lambda^p_n=(\frac {p-L}{p})(\Lambda_n+a\lambda_n)^{p-1}.
\end{equation*}
   We shall leave the details to the reader and we point out that
   as in the case of weighted mean matrices in \cite{G7}, we
   deduce from Theorem \ref{3.1} the following
\begin{cor}
\label{cor1}
   Let $1<p<\infty$ be fixed in \eqref{3.1}. Let $a$ be a constant such that $\Lambda_n+a\lambda_n>0$ for all $n \geq 1$.
  Let $0<L<p$ be a positive constant such that the following
  inequality is satisfied for all $n \geq 1$ (with $\Lambda_0=\lambda_0=0$):
\begin{equation*}
  ( \frac {p-L}{p} )(1+a\frac
  {\lambda_n}{\Lambda_n})^{p-1}+\frac
  {\Lambda_n}{\lambda_{n+1}} \leq \frac
  {\Lambda_n}{\lambda_{n}}(1+a\frac
  {\lambda_n}{\Lambda_n})^{p-1}((1-\frac {L}{p})\frac {\lambda_n}{\Lambda_n}+\frac {\Lambda_{n-1}}{\Lambda_n}+a\frac {\lambda_{n-1}}{\Lambda_n})^{1-p}.
\end{equation*}
   Then inequality \eqref{3.1} holds with
    $U_{p} \leq (p/(p-L))^p$.
\end{cor}

   We now apply the above corollary to the special case of
   \eqref{3.1} with $\lambda_n=\alpha n^{\alpha-1}, \Lambda_n=n^{\alpha}$ for some $\alpha>1$. On taking
   $L=1/\alpha$ and $a=0$ in Corollary \ref{cor1} and setting $y=1/n$, we see that inequality
   \eqref{1.2} holds as long as we can show for $0 \leq y \leq 1$,
\begin{equation}
\label{3.2}
  \Big((1-\frac {1}{p \alpha})\alpha y+(1-y)^{\alpha} \Big
)^{p-1}\Big ( (1-\frac {1}{p \alpha})\alpha y+(1+y)^{1-\alpha}
\Big ) \leq 1.
\end{equation}
   We note first that as $(1-\frac {1}{p \alpha})\alpha
   y+(1-y)^{\alpha} \leq (1-\frac {1}{p \alpha})\alpha
   y+(1+y)^{1-\alpha}$, we need to have $(1-\frac {1}{p \alpha})\alpha
   y+(1-y)^{\alpha} \leq 1$ in order for the above inequality to
   hold. Taking $y=1$ shows that it is necessary to have
   $\alpha \leq 1+1/p$. In particular, we may assume $1< \alpha \leq 2$ from now on and it then follows from Taylor expansion that in order
   for \eqref{3.2} to hold, it suffices to show
\begin{equation}
\label{3.3}
  \Big(1-\frac {1}{p}y+\frac {\alpha(\alpha-1)}{2}y^2 \Big
)^{p-1}\Big ( 1+(1-\frac {1}{p})y+\frac {\alpha(\alpha-1)}{2}y^2
\Big ) \leq 1.
\end{equation}
   We first assume $1<p \leq 2$ and in this case we use
\begin{equation*}
  \Big(1-\frac {1}{p}y+\frac {\alpha(\alpha-1)}{2}y^2 \Big
)^{p-1} \leq 1+(p-1)(-\frac {1}{p}y+\frac
{\alpha(\alpha-1)}{2}y^2)
\end{equation*}
   to see that \eqref{3.3} follows from
\begin{equation*}
  h_{1,\alpha, p}(y):=\frac {\alpha(\alpha-1)p}{2}-(1-1/p)^2+\frac
  {\alpha(\alpha-1)(p-1)}{2p}(p-2)y+\frac
  {\alpha^2(\alpha-1)^2}{4}(p-1)y^2 \leq 0.
\end{equation*}
   We now denote $\alpha_1(p)>1$ as the unique number satisfying $h_{1, \alpha_1,p}(0)=0$ and $\alpha_2(p)>1$ the unique number satisfying
   $h_{1, \alpha_2,p}(1)=0$ and let $\alpha_0(p)=\min (\alpha_1(p), \alpha_2(p)
   )$. It is easy to see that both $\alpha_1(p)$ and $\alpha_2(p)$ are $\leq 1+1/p$ and that for $1<\alpha \leq \alpha_0$, we have $h_{1, \alpha,
   p}(y)\leq 0$ for $0 \leq y \leq 1$.

   Now suppose that $p>2$, we recast \eqref{3.3} as
\begin{equation}
\label{3.4}
   1+(1-\frac {1}{p})y+\frac {\alpha(\alpha-1)}{2}y^2
  \leq \Big(1-\frac {1}{p}y+\frac {\alpha(\alpha-1)}{2}y^2
\Big )^{1-p}.
\end{equation}
    In order for the above inequality to hold for all $0 \leq y
    \leq 1$, we must have $\alpha(\alpha-1)y^2/2 \leq y/p$.
    Therefore, we may from now on assume $\alpha(\alpha-1) \leq
    2/p$. Applying Taylor expansion again, we see that \eqref{3.4}
    follows from the following inequality:
\begin{equation*}
   1+(1-\frac {1}{p})y+\frac {\alpha(\alpha-1)}{2}y^2
  \leq 1+(1-p)(-\frac {1}{p}y+\frac {\alpha(\alpha-1)}{2}y^2)+p(p-1)(-\frac {1}{p}y+\frac {\alpha(\alpha-1)}{2}y^2)^2/2.
\end{equation*}
   We can recast the above inequality as
\begin{equation*}
   h_{2,\alpha, p}(y):=\frac {\alpha(\alpha-1)p}{2}-\frac {1-1/p}{2}+\frac {(p-1)\alpha(\alpha-1)y}{2}-\frac {p(p-1)\alpha^2(\alpha-1)^2y^2}{8} \leq 0.
\end{equation*}
   We now denote $\alpha_0(p)>1$ as the unique number satisfying $\alpha(\alpha-1)
   \leq 2/p$ and
   $h_{2, \alpha_0,p}(1)=0$. It is easy to see that for $1<\alpha \leq \alpha_0$, we have $h_{2, \alpha,
   p}(y)\leq 0$ for $0 \leq y \leq 1$. We now summarize our result
   in the following
\begin{theorem}
\label{thm3.2}
   Let $p>1$ be fixed and let $\alpha_0(p)$ be defined as above, then inequality \eqref{1.2} holds for $1<\alpha
   \leq \alpha_0(p)$.
\end{theorem}

   As we have explained in Section \ref{sec 1}, the study of \eqref{1.2} is motivated by \eqref{1.3}.
   As \eqref{1.2} implies \eqref{1.3} and the constant $(\alpha p/(\alpha p-1))^p$ there is best possible (see
   \cite{G8}), we see that the constant $(\alpha p/(\alpha p-1))^p$ in
\eqref{1.2} is also best possible. More generally, we note that
inequality
   (4.7) in \cite{G8} proposes to determine the best possible constant $U_p(\alpha, \beta)$
   in the following inequality (${\bf a} \in l^p, p>1, \beta
\geq \alpha \geq 1$):
\begin{equation}
\label{3.5}
 \sum^{\infty}_{n=1}\Big{|}\frac
1{\sum^n_{k=1}L^{\alpha-1}_{\beta}(k,
k-1)}\sum^n_{i=1}L^{\alpha-1}_{\beta}(i, i-1)a_i \Big{|}^p \leq
U_p(\alpha, \beta)\sum^{\infty}_{n=1}|a_n|^p.
\end{equation}
    We easily
    deduce from Theorem \ref{thm3.2} the following
\begin{cor}
   Keep the notations in the statement of Theorem \ref{thm3.2}.
   For fixed $p>1$ and $1<\alpha \leq \alpha_0(p)$,
   inequality \eqref{3.5} holds with $U_p(\alpha, \beta)=(\alpha p/(\alpha p-1))^p$ for any $\beta \geq
   \alpha$.
\end{cor}


\noindent {\bf Acknowledgements.}
  The author is supported by a research fellowship from an Academic
Research Fund Tier 1 grant at Nanyang Technological University for
this work.

\end{document}